\documentclass[conference]{IEEEtran}
\IEEEoverridecommandlockouts
\usepackage{amsthm,amsmath,amssymb,mathtools}
\usepackage{graphicx}
\usepackage{cite}

\newtheorem{remark}{Remark}
\usepackage{algorithm}
\usepackage{algorithmic}
\usepackage{booktabs,lipsum}
\usepackage{color}
\DeclareMathOperator*{\argmin}{arg\,min}

\def\BibTeX{{\rm B\kern-.05em{\sc i\kern-.025em b}\kern-.08em
    T\kern-.1667em\lower.7ex\hbox{E}\kern-.125emX}}

\newcommand{\ba}{{\bar a}}

\newcommand{\bT}{{\bar T}}

\newcommand{\br}{{\bar r}}
\renewcommand{\r}{{\bar r}}

\newtheorem{thm}{Theorem}

\newcommand{\ignore}[1]{}

\begin{document}

\title{Continuous-Time Markov Decision Processes with Controlled Observations \\
}

\author{
\IEEEauthorblockN{Yunhan Huang}
\IEEEauthorblockA{\textit{Department of Electrical Engineering} \\
\textit{New York University}\\
Brooklyn, USA \\
yh.huang@nyu.edu}
\and
\IEEEauthorblockN{ Veeraruna Kavitha}
\IEEEauthorblockA{\textit{IEOR Department} \\
\textit{IIT Bombay, Powai}\\
Mumbai, India \\
vkavitha@iitb.ac.in}
\and
\IEEEauthorblockN{ Quanyan Zhu}
\IEEEauthorblockA{\textit{Department of Electrical Engineering} \\
\textit{New York University}\\
Brooklyn, USA \\
quanyan.zhu@nyu.edu}
}

\maketitle
\thispagestyle{plain}
\pagestyle{plain}

\begin{abstract}
In this paper, we study a continuous-time discounted jump Markov decision process with both controlled actions and observations. The observation is only available for a discrete set of time instances. At each time of observation, one has to select an optimal timing for the next observation and a control trajectory for the time interval between two observation points. We provide a theoretical framework that the decision maker can utilize to find the optimal observation epochs and the optimal actions jointly. Two cases are investigated. One is gated queueing systems in which we explicitly characterize the optimal action and the optimal observation where the optimal observation is shown to be independent of the state. Another is the inventory control problem with Poisson arrival process in which we obtain numerically the optimal action and observation. The results show that it is optimal to observe more frequently at a region of states where the optimal action adapts constantly. 
\end{abstract}


\section{ Introduction }
Markov Decision Processes (MDP) \cite{Puter} are widely applicable to many real world problems including queueing systems \cite{stidham1985optimal}, communication systems \cite{altman2002applications}, and motion planning of robotics \cite{ragi2013uav} etc. Many MDP frameworks assume that continuous updates of direct or indirect observation is available. This assumption becomes inadequate to capture many applications where observations are either limited or are costly. For example, in many cyber-physical systems, controllers and sensors are physically separate. Control needs to be incessantly applied based on the remote sensing of the physical system. However, sensors cannot always provide continuous sensing due to limited battery capacity or weak processor. As another example  for large-scale networks, e.g., Internet of Things (IoT) networks, continuous sensing is costly and unnecessary.

To address these issues, we propose an MDP framework with controlled and limited observations. We consider in this paper a continuous-time jump MDP with discounted cost criteria that take into account the cost of observations. In the framework, the decision maker cannot observe the state continuously but has to determine the next observation time after one observation and an action/control trajectory during the time between two observation points. Hence, the decision maker has to jointly determine the control trajectory and observation points. There is a fundamental trade-off between actions and observations. On one hand, more information/observations facilitate the decision maker to design better action so as to increase the system performance. On the other hand, more information/observation requires a higher cost of communication and power usage.

In this paper, we establish a theoretical framework that the decision maker can utilize to find the optimal observation epochs and the optimal actions jointly. To be more specific, we use dynamic programming techniques to characterize the value function and leverage value iterations to compute it numerically. The dynamic programming equation is obtained by forming a consolidated utility between two observation points. At each iteration, one has to solve an optimal control problem (or a finite-dimensional optimization problem, depending on the structure of the action) whose terminal cost consists of the transition cost at the next observation and the cost of observation. 

We provide two case studies to illustrate our framework: one is the gated queueing system for which we obtain an explicit dynamic programming equation and characterize explicitly both the optimal observation and the optimal action. The other one is an inventory control problem with Poisson arrival and departure processes where we perform value iterations and obtain numerically the optimal observation points and the optimal actions. The results show that at a region of states where the optimal action adapts constantly, it is optimal for the decision maker to observe more frequently. At a region of states where optimal action barely changes, the decision maker observes less frequently.

\textit{Related work:} One related work is Partially Observable Markov Decision Process \cite{krishnamurthy2016partially} (POMDP). POMDP deals with partial observations with the state. States are continuously observed in an indirect fashion. However, this work deals with the scenario where observations are only available at a discrete set of time instances over a continuous time-line. And the decision maker has to decide the set of time instances for observation. Another related work is networked control systems where sensing and control signals are communicated through networks. This framework is also related to optimal control of dynamical systems under intermittent measurements \cite{bacsar1995minimax,ades2000stochastic,imer2006optimal}. But these works have focused on linear dynamical systems whose dynamics are governed by differential/difference equations and the control aspect of the problem, i.e., finding controllers under a given pattern of observations. For example, in \cite{ades2000stochastic}, the authors aim to find an optimal control for time-varying linear dynamical systems when observation instants are Poisson distributed. In this work, we aim not only to find optimal control policies but also to characterize jointly the optimal sensing policies under costly measurements. To the best of our knowledge, this is the first work studying Markov decision process with both controlled action and observation.

\textit{Organization of the paper:} In Section \ref{SecProDes}, the problem of continuous-time jump MDP with controlled observation is formulated. We analyze the problem and develop a general theoretic framework in Section \ref{SecAnaysis}. A gated polling system is presented in Section \ref{sec_queue} where its optimal solution is characterized explicitly. In Section \ref{sec_example}, we study an inventory control problem under the framework whose optimal observation time is obtained numerically.
 
\section{Problem Description}\label{SecProDes}
Consider  a controlled Markov Jump process (MJP)\footnote{
Markov Jump Process (MJP) is specified by the rates of the exponential transition times for the given state, $\lambda_x $  for all $x \in S$ the (countable) state space, and the transition probabilities $q(x' |  x, a)$ at the time of jump.  The (exponentially distributed)  transition times and the state transitions are independent.
}   $\{X(t) \}$,    whose transitions  can be controlled.  Often it is expensive to observe  and control the process frequently and  this is the main theme of this paper.   
The controller  further wants to control the observation epochs, which also become the  decision epochs.

{\bf Controls:}    Let $A(t)$ represent the control process which determines the rate of transitions and or the transition probabilities.
For example, this can represent control  (only) of  rate of transitions;
given that the MJP  is in state $x$, its next transition epoch is given by a (time-varying rate) Poisson process with rate process $\lambda_x a(\cdot)$, when the rate control is the measurable function $a (\cdot) \in L^\infty [0, \infty ]$.  Here, $\{\lambda_x \}_{x \in S}$ are fixed.  
Alternatively, one can control the transition probabilities.  
 We can also have observation rate control, specified by $T_x$  which specifies the time duration till the next observation epoch when the process is in state $x$.

{\bf Policy:}  We consider stationary Markov policies; for example, a policy can be of the following form (for some appropriate  $p$):  $$
\pi = \left \{ \left ( a_x (\cdot),  T_x  \right )  \in   L^p [0, \infty)  \times [\underline{T}, \infty) :  \mbox{ for any }   x \in S \right  \},
$$where $S$ is the state space, $\underline{T}$ is the minimum time gap between two  successive observations and $L^p$ is the space of $p$-integrable functions.  Basically any policy prescribes an action to be chosen at every observation epoch, based on the state of the system at that time point, and the time till the next observation epoch. 
The action could  be a simple single decision, i.e.,   $a(t) = a$ for all $t$ for some action $a$ (as is usually the case with Markov decision process), or it could be a (most general)   open loop control (a measurable function of time)  which is applicable till the next information update.  
 We  show that these (Markovian) policies are sufficient  under the assumptions of this paper. 

{\bf Transition probabilities:} Let 
\begin{eqnarray*}
q(x, x' ; a, T) := \hspace{-18mm} &\\ &
Prob(X(T) = x' |  X (0)  = x,  A( [0, T]) = a(\cdot), T_1 = T)
\end{eqnarray*} represent the probability that the MJP  is in state  $x' \in S$ after time $T$, given that the initial condition is $X (0)  =x$ and  given that the open loop control $a(\cdot)$ is chosen   till the second observation epoch which is after $T$ units of time.  Here $A([a_1, a_2])$ represents the control for the time duration $[a_1, a_2].$

{\bf Utilities:}  As is the case with Markov decision processes, the overall utility is that combined over several time periods.  However,  we have insufficient/infrequent observation epochs and a state-blind decision has to be chosen for any time period spanning between two observation points. Thus it is appropriate to consider consolidated utilities,  consolidated for time periods between subsequent observation epochs.  
Let $\bT_k := \sum_{i < k}  T_i$ (and $\bT_1 = 0$) represent the instance of $k$-th observation, where $T_k$ is the time gap between $k$-th and $(k+1)$-th observation epochs. 
Let $$\br_k := \br  (X(\bar{T}_k), A(\bar{T}_k + \cdot), T_k)$$ represent the consolidated utility that corresponds to the time period between $k$-th and  $(k+1)$-th observations epochs. 
Below we describe few examples.

\underline {Example:} Let $r(X(t), A(t))$ represent the instantaneous utility at time $t$    and then the discounted utility  is given by:
\begin{eqnarray}
J(\pi, x) = E^{\pi}_{x} \left [ \int_0^\infty \beta^t  r(X(t), A(t) ) dt  + \sum_i  \beta^{ \bT_i } g (T_i ) \right ]. 
\label{Eqn_JFunction}
\end{eqnarray}
where the expectation is with respect  to the given policy  $\pi$ and  initial condition  $X(0) = x$ and where $g(T_i)$ is the cost for choosing observation epoch $T_i$ which is supposed to be bounded and decreasing in a admissible region of $T_i$. 
 Assume $|r (\cdot, \cdot) | < B$, i.e., uniformly bounded  utilities for discounted cost examples.  
In this case one can define the consolidated utilities as  below:

\vspace{-3mm}
{\small \begin{eqnarray}
\br_k &=& \br (X(\bar{T}_k), A(\bar{T}_k + \cdot), T_k)    \label{Eqn_Integral_utils}\\
 & =&  E_x^\pi\bigg  [ \int_0^{T_k}  \beta^t  r \left (X(\bar{T}_k+t), A(\bar{T}_k+t)  \right ) dt \nonumber 
\\ 
&&  \hspace{40mm}  \bigg  |  X(\bar{T}_k),  A(\bar{T}_k +  \cdot ), T_k  \bigg ].  \nonumber
\end{eqnarray}
By Fubini's Theorem and Markov property \cite{durrett2019probability}:
\begin{eqnarray}
\br (x, a(\cdot), T) \hspace{-18mm}  \label{Eqn_Fubini}  &\\ & = 
  \int_{0}^{T} \beta^t E^{a (\cdot) }_{x} \left [   r(X(t), A(t) ) \bigg |  X(0 ) = x,   A(  \cdot ) = a(\cdot), T    \right ]   dt .  \nonumber
\end{eqnarray}}
The above is the time integral of the expected utilities of the Markov  jump process  till the (chosen) observation epoch $T$, when started with $x$ and when open loop control $a(\cdot)$ is used till $T$.

\underline{Other examples:}
One can alternatively consider discrete time average/sum  of a measurement of the system over an observation period as the consolidated utility. 
In all, the consolidated utility is either discrete/continuous average/sum of the measurements of the system during the corresponding observation period. 
The consolidated utilities can depend upon the dynamics of the  system over an inter-observation period   in various  other ways.   One may be interested in the probability that the system crosses a certain threshold, or one might be interested  in extreme values of a measurement related to the system during the said period.  In the queueing system based example, considered in section \ref{sec_queue}, the consolidated utilities correspond to the sum of the waiting times of the customers that arrived in one observation period. 

\section{Problem Formulation and analysis}\label{SecAnaysis}
To begin with, 
we are interested in optimizing the following accumulated utility, constructed using   discounted   values of several consolidated utilities:
\begin{eqnarray}
v(x) &=& \sup_\pi J(\pi, x) \mbox{ with } \label{Eqn_J_pi}  \\
J(\pi, x)  
&=&
\sum_{k=1}^\infty  E^{\pi}_{x} \left [ \beta^{\bT_k} \bigg (  \br_k  + g (T_k ) \bigg  )   \right ].   \nonumber
\end{eqnarray}
Observe here that consolidated  utilities   given by  (\ref{Eqn_Integral_utils}) when accumulated as above would generate  the utility  exactly as in (\ref{Eqn_JFunction}).

\ignore{
{\bf Analysis:}  One can rewrite the above cost as in the following by using the tower property and by Markovian property:
\begin{eqnarray}
J(\pi, x) & =&   \sum_i E^{\pi}_{x} \left [ \beta^{\bT_i}  \left ( \int_0^{T_i} \beta^t    r(X(t+\bT_i), A(t+\bT_i) ) dt + g (T_i )  \right )  \right ]  \nonumber  \\
&=&
\sum_i E^{\pi}_{x} \left [ \beta^{\bT_i} \bigg (  \r  ( X(\bT_i),  \   A(\bT_i + \cdot ), \  T_i  ) + g (T_i ) \bigg  )   \right ]  \mbox { with} \label{Eqn_J_pi}  \\
\r (x, a (\cdot), T ) & := & E^{\pi}_{x} \left [    \int_{0}^{T} \beta^t   r(X(t), A(t) ) dt  \bigg |  X(0 ) = x,   A(  \cdot ) = a(\cdot)    \right ]  \nonumber \\
& = &    \int_{0}^{T} \beta^t E^{a (\cdot) }_{x} \left [   r(X(t), A(t) ) \bigg |  X(0 ) = x,   A(  \cdot ) = a(\cdot)    \right ]   dt ,
\label{Eqn_Discrete_JFunction}
\end{eqnarray}by Fubini's theorem.}
A close look at  (\ref{Eqn_J_pi})  shows that this is   a discounted cost discrete-time  MDP, with discount factor $\underline{\beta} := \beta^{\underline{T}}$,   Markov state and Markovian actions given respectively by
$$
Z_k :=  (X(\bT_k),  {\tilde T}_k ),  A_k = (a_k (\cdot), T_k) \mbox{ with }   {\tilde T}_k:= \bT_k-(k-1) \underline{T}
$$  
and running cost  equal to
$$
r(Z_k, A_k) =  \beta^{ {\tilde T}_k}      \bigg (  \r   \big ( X(\bT_k),  \   a (\cdot), \  T_k \big  ) + g (T_k ) \bigg  ) . 
$$ 
That is,  utility in (\ref{Eqn_J_pi}) is given by
\begin{eqnarray*}
J(\pi, x) & =& \sum_{k= 1}^\infty  \underline{\beta}^{k-1}   r(Z_k, A_k),  \mbox{ with } \underline{\beta} := \beta^{\underline{T}}. 
\end{eqnarray*}
Thus one can  derive  relevant results from the current MDP literature, mainly the results available to Polish (Banach and separable) spaces {(e.g., \cite{Puter}).} For example, we can conclude the sufficiency of stationary Markov policies. 
The value iteration converges for any $\beta <1$.
We further observe that the optimal policies need not depend upon the observation instances, $\{\bT_k\}$ in the following:

\begin{thm}[{\bf Dynamic programming equation}]
\label{Thm_DP_eqns}
The value function defined by (\ref{Eqn_J_pi}) satisfies the following dynamic programming equation:
\begin{eqnarray*}
v(x) &=&   \sup_{ a \in L^p  [0, \infty), T  \in [{\underline T}, \infty) }
  \bigg \{ \br (x, a(\cdot), T)   \\  && \hspace{5mm}+ 
  \beta^T     \sum_{x' \in S}  q(x, x';   a, T)  v(  x' )  + g(T) \bigg  \}.
\end{eqnarray*}
If there exists a policy $\pi^* = \{ (a^*_x (\cdot ),  T^*_x ) : x \in  S \} $ such that
\begin{eqnarray}
\label{Eqn_DP_eqns}
v(x) \hspace{-1mm}&\hspace{-1mm}=\hspace{-1mm}& \hspace{-1mm}
  \bigg\{    \br (x, a_x^*(\cdot), T_x^*)  
  \nonumber \\ &&
 +  \beta^{T_x^*}     \sum_{x'  \in S}  q(x, x';   a_x^*,  T_x^*)   v(  x' )  + g(T_x^*) \bigg  \},     \hspace{8mm}
\end{eqnarray}for all $x \in S$,  then  $\pi^*$ is the optimal policy.
\end{thm}
{\bf Proof: }  It immediately follows from the discussion above the theorem and by noting that the actual  value function (after definitions) equals
$$ 
v(\bT, x) = \beta^{\bT}  v(x).   \hspace{30mm}  \blacksquare
$$ 

{\bf Remarks:}
 1) One may want to control only $\{T_x\}$ and not $\{a_x(.)\}$, or vice versa or the both.
 
 2) For examples as in (\ref{Eqn_JFunction}),
in view of (\ref{Eqn_Fubini}),     we need to compute  the quantities   $  E_x   \left [  r(X(\tau), a(\tau ) )  \right ]$,  for all $\tau$ within any given observation period.
This can be computed if one  can estimate the expectation of the MJP at any given time $\tau$ which started with $X(0)  = x$  (for some $x$) and whose transitions are governed by $a(\cdot).$

3) One can solve these   dynamic programming equations  utilizing the usual  value iteration method;  Given the $k$-th estimate of  value function,  $\{ v_k (x)\}_{x \in S}$, 
use (\ref{Eqn_DP_eqns}) to get the next estimate $\{ v_{k+1} (x)\}_{x \in S}$ and repeat this till the fixed points   converge. This iterative method is  guaranteed to converge when $\beta < 1$, in view of Theorem \ref{Thm_DP_eqns} and the discussions prior to that.  More details are provided in the immediate subsection.

\subsection{Solving DP Equations via Value Iterations}
When one considers examples defined through integrals as in (\ref{Eqn_JFunction}),  in view of (\ref{Eqn_Fubini}),   the  DP equations (for any given $t$)  can be solved by solving an appropriate optimal control problem as explained below.  For such examples, the DP equations can be rewritten as the following:

\vspace{-4mm}
{\small \begin{eqnarray*}
\sup_{a (\cdot), T }
\Bigg \{  
  \int_{0}^{T} \beta^t  c(t;  a ) dt    
   +  \beta^{T}     \sum_{x'  \in S}  q(x, x';   a,  T)   v(  x' )  + g(T)  \Bigg \} ,  \\
  c(t; a) :=  E^{a (\cdot) }_{x} \left [   r(X(t), A(t) ) \bigg |  X(0 ) = x,   A(  \cdot ) = a(\cdot)    \right ]   dt.
\end{eqnarray*}} 
The value iteration for such a case can be implemented using the following iterative procedure:
\begin{eqnarray*}
v_{k+1} (x ) = 
\sup_{a (\cdot) }
\Bigg \{  
  \int_{0}^{T} \beta^t    c(t; a)  dt  \hspace{40mm}  \\
   +  \beta^{T}     \sum_{x'  \in S}  q(x, x';   a,  T)   v_k (  x' )  + g(T)  \Bigg \}, 
\end{eqnarray*}
after appropriately initializing $\{v_0(x)\}$ and the iteration is stopped at some $k$ after testing for proper convergence. 
It is easy to observe that each stage of the above iterative procedure and for each $T$, can be solved using optimal control tools and this procedure is explained using an example  
in Section \ref{sec_example}.

\subsection{Separable Utilities}
If the  consolidated  utilities, consolidated  over an observation period,  are  separable as below:
\begin{eqnarray}
\label{Eqn_sep_cost}
\br (x, a(.), T) = \br_a (x, a(.) )  +  \br_T (x, T). 
\end{eqnarray}
And further say the transition function does not depend upon the control $a(\cdot)$. 
In this case,       the DP equations simplify to the following:
\begin{eqnarray*}
v(x) =   \sup_{ a \in L^p [0, \infty), T  \in [{\underline T}, \infty) }
  \bigg  \{      \r_{x,a}  ( x,   a( \cdot )   ) +  \r_{T}  (x, \ T)    \\   \hspace{5mm}+
  \beta^T     \sum_{x'}  q(x, x';  T)    v(  x' )  + g(T) \bigg  \}.  \end{eqnarray*}
  This is equivalent to the following: 
  
\vspace{-4mm}  
  {\small \begin{eqnarray*}
  \sup_{   T  \in [{\underline T}, \infty) }
  \left \{      \r^*_{x}  ( x  ) +  \r_{T}  (x, \ T)    +
  \beta^T     \sum_{x'}  q(x, x'; T)  v(  x' )  + g(T) \right \},  \\
  \end{eqnarray*}}
  with
\begin{eqnarray*}
  \r^*_x  (x)  :=    \sup_{ a \in L^p [0, \infty) }   \r_{x,a}  ( x,  \  a( \cdot )   ),  \hspace{30mm}
\end{eqnarray*}which is like DP equations of the usual  MDP  with control as $T$.

\section{ Case Study:  Gated Queueing Systems}
\label{sec_queue}

In this section, we  consider an example that has separable utilities.  Further the consolidated utilities of this example are  not of integral form as in (\ref{Eqn_JFunction}).

Consider a queuing system which operates like a gated polling system where the gate openings are controlled dynamically. We further wish to control the server speed based on the number of customers. Here, the observation epochs are the gate opening instances. One can think of these systems to have two waiting rooms;  when the gate is opened all the customers from outer room enter the inner room and the gate is closed immediately after.  A server speed is allocated to serve the workload that entered the inner room. In addition the next gate opening instance  is decided and  the arriving customers accumulate in the outer room till the gate opens again. 

The arrivals to this system are modelled by a Poisson process with constant rate $\lambda.$  The service times are independent and identically distributed across customers (when they  are served at the same speed);  however,  the expected  service  time is inversely proportional to the server speed.

We are interested in optimizing a discounted cost related to the waiting times of the customers.  Any customer has to wait (in outer room)  for the next instance at which the gate  opens, and then  will have to wait for its turn in the inner waiting time.  We assume that the customers are served using First come first serve discipline. 
The discount factor is the same for all the customers waiting  in an observation period, i.e., we are interested in optimizing the expected value of the  following:
$$
\sum_{k=1}^\infty \beta^{\bT_k}  E[ W_k ],  
$$where  $W_k$ is the sum total of the waiting times of all the customers waiting during  $k$-th observation  period, i.e., during the time period between $(k-1)$-th and $k$-th gate opening instances. 
These waiting times are due to  the customers waiting in the outer  and  inner rooms.  The  customers waiting in the outer room are the  new arrivals into the system; the number of  new arrivals ($N_k$) is Poisson distributed with parameter $\lambda T_k$, if  $T_k$ is the  length of the $k$-th observation period. 
The expected waiting time of all the customers waiting in outer room during this period\footnote{The  residual time till the next epoch (i.e., the waiting time till the gate opens) of an   un-ordered arrival  that arrived during an interval of length $t$, is uniformly distributed over $t$ and the expected number of arrivals equals $\lambda t$. } thus equals:
$$
\lambda T_k^2 / 2 .
$$
If the server speed is $a_k$   and  if $X_{k-1}$ (observe $X_{k-1} = N_{k-1}$)  number of customers   enter  the inner room at $(k-1)$-th gate open instance, then the sum of expected waiting times of these customers equals:
$$
  \frac{1}{a_k}  +  \frac{2}{a_k}  +  \cdots + \frac{3}{a_k} + \cdots   + \frac{X_{k-1}}{a_k}.
$$Basically the $i$-th served customer has to wait for the service of  $i$-customers, before departing. Thus
$$
E[W_k | X_{k-1}, T_k, A_k = a_k] = \frac {\lambda T_k^2 } {2 } +  \frac{ X_{k-1}^2 + X_{k-1}}{2 a_k}  .
$$
 Furthermore,  there is a cost for (frequent) observations   $g(T_k)$  and a cost for server speed $\eta(A_k)$.     We make a simplifying assumption that the load factor is moderate and that the customers in the inner room are served with probability close to one  before the minimum observation period  $\underline{T}.$
We would like to optimize the following 
$$
\sum_{k=1}^\infty E^\pi_{0}\left [  \beta^{ \bT_k }   \bigg  ( \frac {\lambda T_k^2 } {2 } +  \frac{ X_{k-1}^2 + X_{k-1}}{2 A_k} + g(T_k) + \eta (A_k)  \bigg ) \right ]
  $$where  $\pi$ is any   stationary Markov policy that specifies the next observation period $T_k = T(x)$ and the server speed $A_k = A(x)$, given  that the number of customers who have entered the inner room at $(k-1)$-gate open instance,  
  $X_{k-1} =x$.

This is an example of separable utilities as in (\ref{Eqn_sep_cost}), with 
$$\br_{a}(x, a) = \frac {x^2+x}{2a} + \eta(a) \mbox{  and   } \br_T (x, T) = \frac {\lambda T^2}{2} .$$ 
Further, the transition probabilities are independent of action $A_k$:
$$
q(x, x'; a, T) = \exp (- \lambda T) \frac{ (\lambda T)^{x'}}{(x')!}.
$$
Assume a linear server speed cost $\eta(a) =\eta a$ and we immediately arrive at  
$$
r^* _{x} (x) =   \sqrt{  2\eta x(x+1)   }     \mbox{ and }      a^*(x) =  \sqrt{  \frac{x(x+1)}{2 \eta }  }.
$$
Finally, the DP equations for this example are given by the following and simplified as below: 
\begin{eqnarray*}
v(x) &=&
    \inf_{   T  \in [{\underline T}, \infty) }
  \bigg  \{   \sqrt{  2\eta x(x+1)   }    +  \frac {\lambda T^2 }{2 }   \\  && \hspace{15mm}  +
  \beta^T     \sum_{k=0}^\infty   e^{- \lambda T}  \frac{ \lambda^k T^k }{k!}   v(  k )  + g( T) \bigg   \}    \\
  &=& \sqrt{  2\eta x(x+1)   }   + 
    \inf_{   T  \in [{\underline T}, \infty) }
  \bigg  \{     \frac {\lambda T^2 }{2 }   + \\  && \hspace{25mm} +
  \beta^T     \sum_{k=0}^\infty   e^{- \lambda T}  \frac{ \lambda^k T^k }{k!}   v(  k )  +g( T)  \bigg    \}.  
  \end{eqnarray*}
 For further analysis of the optimal policy, one needs to solve the remaining optimization problem using numerical methods. But it is interesting to observe that the 
 server speed depends upon the number of customers who have  entered the inner room, while the optimal observation epochs are independent of the state of the system, i.e.,
 $$
 A^*_x = \sqrt{  \frac{x(x+1)}{2 \eta }  } \mbox{ and }  T^*_x =   T^*,  \mbox{ for all } x,
 $$for some constant $T^*$.
 
One can consider interesting variants of this problem. For example, one can consider open-loop control for the server speed and choose a sequence of server speeds to be used for each of the $X_k$ waiting customers. One can consider a constraint on the next observation period, which has to be longer than the time taken to complete all the previous jobs, in which case the structure of separable utilities disappears.

\section{Case Study: Inventory Control}
\label{sec_example}

\subsection{Setup of the case study}

In this case study, we consider an energy harvesting type of inventory control problems with the aim to maintain an inventory level close to $\theta$ where both the arrivals and the departures are considered random and every arrival/departure corresponds to one unit. The departures are coming  according to a Poisson process with fixed rate $\mu$, while the arrivals are modelled by inhomogeneous Poisson process with a controlled time-varying rate $a(\cdot)$. 

{At each observation epoch, we have a fixed time interval $[0,T]$.} Let $X_k$, be the amount of inventory at the end of $k$-th observation. With every arrival, one unit is added to the inventory and one unit is reduced with every demand departure. Thus, the amount of inventory at any time point before $(k+1)$-th epoch is given by
$$
X_k(t) = X_k + \mathcal{A}(t;a) - \mathcal{D}(t)=X_k + \sum_{i=1}^{\mathcal{N}(t;a)} \xi_i
$$
where $\mathcal{A}(t;a)$ and $\mathcal{D}(t)$ are the number of arrivals and departures by time $t$, respectively, since the last observation and $\mathcal{N}(t;a)$ is the total number of arrivals and departures with $\xi_i=1$ if it is arrival and $-1$ if it is departure. Further, we assume that the departures can be backlogged infinitely, if the inventory is empty. Hence, we have $X_k\in\mathbb{Z}$ and $X_k(t)\in\mathbb{Z}$ for all $t$.

The overall utility depends on the cost spent on the acceleration process and the deviation from the targeted inventory $\theta$ which is
$$
\int_0^\infty \beta^t \left\{ E\left[(X_k(t)-\theta)^2 \right]  + \nu a(t) \right\}dt
$$
where $\nu\geq 0$ characterizes the cost of accelerating the arrivals. 

\subsection{Analysis}
Note that 
$$
E[\xi_i] =   \frac{   \ba(t) - \mu t } {\ba(t) + \mu t} ,  \mbox{ with } \ba(t) = \int_0^t a(s) ds,
$$
where $\bar{a}(0)=0$ and $a(t)\in[0,\bar{\ba}]$ for every $t$. Then by Wald's lemma \cite{wald1945some}, the deviation cost can be rewritten as 
{\small
$$
\begin{aligned}
& E\left[\left ( X_k (t) - \theta  \right )^2 \right]\\
=& ( X_k - \theta )^2  + 2 (X_k - \theta)   E[{\cal N} (t; a)]  \frac{   \ba(t) - \mu t } {\ba(t) + \mu t} + E \left[ \left(\sum_{i=1}^{{\cal N} (t; a)}  \xi_i  \right)^2  \right]\\
=&  ( X_k - \theta )^2  + 2 (X_k - \theta)     (\ba(t) + \mu t)  \frac{   \ba(t) - \mu t } {\ba(t) + \mu t}\\
&+ \left [     E[{\cal N} (t; a)]  E[\xi^2] +    E[{\cal N} (t; a) (  {\cal N} (t; a)  - 1) ]   \left ( E[\xi_i] \right )^2 \right ]\\
=&  ( X_k - \theta )^2  +  2 (X_k - \theta)     (\ba(t) + \mu t)  \frac{   \ba(t) - \mu t } {\ba(t) + \mu t}\\
&+ \left [     E[{\cal N} (t; a)]\cdot 1 +   (\ba(t) + \mu t)^2 \frac{   (\ba(t) - \mu t )^2} {(\ba(t) + \mu t)^2} \right ]  \\
=&  ( X_k - \theta+( \ba(t) - \mu t  ) ) ^2   + (\ba(t) + \mu t).\\
\end{aligned}
$$}

By Theorem $1$, the DP equation is
\begin{equation}\label{DPEqExamp}
\begin{aligned}
v(x) =& \inf_{a\in L^\infty[0.\infty),T} \Bigg\{  \int_{0}^T  \beta^t r(x,a(t),t)   dt   &\\
&+  \beta^T  \sum_{x'}q(x, x’; a, T)v(x')+g(T)   \Bigg\},
\end{aligned}
\end{equation}
where
$$
r(x,a(t),t)=( x - \theta+( \ba(t) - \mu t  ) ) ^2   + (\ba(t) + \mu t)-a(t)\nu,
$$
and the transition probability given control $a(\cdot)$ and time period $T$ would equal
\begin{eqnarray*}
\begin{aligned}
&\hspace{-1mm}q(x, x'; a, T) = \textrm{Prob} (  X_k(T)  = x'   | X_k = x, a (\cdot) )\\
& \hspace{5mm}= \left \{
\begin{array}{llll}
\sum_{k= x'-x}^\infty  \frac{ e^{-\ba(T) } (\ba(T))^k }{k!}  \frac{ e^{-\mu T} (\mu T)^{k-x'+x}  }{(k-x'+x)!}  &\mbox{ if }  x' > x \\
\sum_{k= x-x'}^\infty  \frac{ e^{- \mu T } ( \mu T)^k }{k!}  \frac{ e^{-\ba(T)} (\ba(T))^{k-x+x'}  }{(k-x+x')!}  &\mbox{ else.}
\end{array} \right .\\
\end{aligned}
\end{eqnarray*}
Here, $g(T)\coloneqq -\kappa T$ characterizes the observation cost with $\kappa>0$.

{\bf Value Iteration:} 
In each iteration (say at iteration $k$), we start with an estimation of the value function $\{v_{k}(x)\}$ and we obtain the new estimates $\{v_{k+1}(x)\}$ by evaluating the fixed point equation (\ref{DPEqExamp}). It is easy to observe that, given $\{v_k(x)\}$, (\ref{DPEqExamp}) can be evaluated by solving the following optimal control problem for each $T$ and choosing the optimal $T$:
\begin{equation}\label{InnerOCExamp}
    v_{k+1}(x)=\inf_{a\in L^\infty[0,T]} \int_0^T \beta^t r(x,a(t),t)dt + h(\bar{a}(T),T),
\end{equation}
where the terminal cost is 
$$
h(\bar{a}(T),T)\coloneqq \beta^T  \sum_{x'}q(x, x'; a, T)+g(T).
$$

\begin{remark}
Note that $x$ is not a state in optimal control problem (\ref{InnerOCExamp}). Let $y(t)\coloneqq \bar{a}(t)$ and treat $y$ as the fictitious state of the optimal control problem. Thus, we have $r(x,a(t),t)= \beta^t (x-\theta +y(t)-\mu t )^2 + (y(t)+\mu t)+a(t)\nu$ with $\dot{y}(t)=a(t)$ and $y(0)=0$.
\end{remark}

Then, at iteration $k$, we need to solve the following optimal control problem (for any given $x$)
{
\begin{equation}\label{OPCon}
\begin{aligned}
\inf_{a\in L^\infty[0,T]} &\int_0^T \beta^t \{(x-\theta +y(t)-\mu t)^2\\
&+(y(t)+\mu t)+a(t)\nu\}dt+ h(y(T),T)\\
s.t.\ \ \ \ &\dot{y}(t)=a(t),\ y(0)=0.
\end{aligned}
\end{equation}}
The Hamiltonian of problem (\ref{OPCon}) is given by
$$
H(t,a,y,\lambda)= \beta^t \{(x-\theta + y - \mu t)^2 + (y+\mu t)+a\nu \} + \lambda a
$$
where $\lambda$ is the costate. With a simple application of minimum principle \cite{liberzon2011calculus}, we know that the optimal solutions $a^\star$ and the corresponding state $y^\star$ need to satisfy the following conditions:
{\small
\begin{equation}\label{NecCon}
\begin{cases}
&\dot{y}^\star(t)=a^\star(t),\\
& \dot{\lambda} = - \beta^t \{ 2(x-\theta +y^\star(t) -\mu t) +1 \},\ \lambda(T)=\frac{\partial h}{\partial y} (y(T),T),\\
& a^\star(t)=\begin{cases}
0\ \ \mbox{if } \beta^t \nu + \lambda >0,\\
\bar{\ba}\ \mbox{otherwise}.\\
\end{cases}\\
\end{cases}
\end{equation}}
\begin{remark}
In the upper level, we have a dynamic programming equation (\ref{DPEqExamp}) which characterizes the value at $x$. In the lower level, we have an optimal control problem whose open-loop solution generates $a_x$ for every $x$ and $T$. The necessary conditions indicate that the optimal control is a bang-bang type of control. At each value iteration, one has to solve problem (\ref{OPCon}) for every $x$. To solve (\ref{OPCon}), one may resort to numerical methods such as nonlinear programming based direct method and indirect methods \cite{Diehl11numericaloptimal}.
\end{remark}

When the Poisson arrival process is homogeneous, i.e., $a(t)$ is fixed over $t$, problem (\ref{InnerOCExamp}) becomes a finite-dimensional optimization problem with 
$$
v_{k+1}(x) = \inf_{a\in[0,\bar{\ba}],T} \int_0^T \beta^t r(x,at,t)dt+h(aT,T),
$$
Thus, we have
{
\begin{equation}\label{VIFixedRate}
\begin{aligned}
&\hspace{-5mm}v_{k+1}(x)\\
=&\min_{a\in[0,\bar{\ba}],T}\int_0^T \beta^t \{ (x-\theta+at -\mu)^2 +at +\mu t\}dt\\
&\ \ \ \ \ \ \ \ \ \ \ \ \ \ +h(aT,T)\\
=&\min_{a\in[0,\bar{\ba}],T} F(a,T;x,v_k), \\
\end{aligned}
\end{equation}
where
\begin{equation}
\begin{aligned}
&\hspace{-4mm}F(a,T;x,v_k)\\
=&a^2\big\{K_2(T)-K_2(0) \}\\
&+a\{[2(x-\theta)+1](K_1(T)-K_1(0))-2\mu(K_2(T)-K_2(0))\\
&+\nu(K_0(T)-K_0(t)) \big\}\\
&+(x-\theta)^2(K_0(T)-K_0(0))\\
&+\mu[1-2(x-\theta)](K_1(T)-K_1(0))+\mu^2(K_2(T)-K_2(0))\\
&+ \beta^T \Big[\sum_{x'=x+1}^\infty v_k(x')    \sum_{k=x'-x}^{\infty}        \frac{e^{-aT}(aT)^k}   {k!}   \frac{e^{-\mu T}(\mu T)^{k-x'+x}}    {(k-x'+x)!}\\
&+\sum_{-\infty}^{x'=x} v_k(x')   \sum_{k=x-x'}^\infty    \frac{e^{-\mu T}(\mu T)^k}   {k!}    \frac{e^{-aT}(aT)^{k-x+x'}}   {(k-x+x')!}\Big] -\kappa T,
\end{aligned}
\end{equation}}
with 
$$
\begin{aligned}
K_0(t)&=\frac{\beta^t}{\ln{\beta}}\\
K_1(t)&=\frac{t\beta^t}{\ln{\beta}}-\frac{1}{\ln{\beta}}K_0(t)\\
K_2(t)&=\frac{t^2\beta^t}{\ln{\beta}}-\frac{2}{\ln{\beta}}K_1(t).
\end{aligned}
$$
The minimization problem admits a minimum since the problem has a continuous objective function over a compact set $[0,\bar{\ba}]\times [\underline{T},\overline{T}]\subset \mathbb{R}^2$. To learn the optimal control $a^\star_x$ and the optimal time for next observation $T^\star_x$, given current observation $X_k=x$, we resort to value iterations given by $(\ref{VIFixedRate})$. Algorithm \ref{ObConAlgo} describes the steps to find $a^*_x$ and $T^*_x$.

\begin{algorithm}
 \caption{VI Based Optimal Observation \& Control Algorithm}
 \label{ObConAlgo}
 \begin{algorithmic}[1]
 \renewcommand{\algorithmicrequire}{\textbf{Input:}}
 \renewcommand{\algorithmicensure}{\textbf{Output:}}
 \REQUIRE Discount factor $\beta$, Departure rate $\mu$, Acceleration cost $\nu$, Reference inventory $\theta$, Tolerance $\epsilon$;
 \ENSURE  Optimal Control $a^\star_x$, Optimal Observation $T_x^\star$, for every $x\in\mathbb{Z}$;
 \\ \textit{Initialization}: $v_0(x)\coloneqq |x-\theta|$, for every $x\in\mathbb{Z}$; $k=0$;
 \WHILE{$k \geq 0$}
  \FOR {$x\in\mathbb{Z}$}
  \STATE $v_{k+1}(x)=\min_{a\in[0,\bar{\ba}],T} F(a,T;x,v_k)$;
  \ENDFOR
  \IF {$\Vert v_{k+1}-v_k \Vert\leq \epsilon$}
  \STATE \textbf{break};
  \ENDIF
  \STATE $k=k+1$;
  \ENDWHILE 
  \FOR{$x\in\mathbb{Z}$}
  \STATE $v^\star(x)=v_{k+1}(x)$;
  \ENDFOR
  \FOR{$x\in\mathbb{Z}$}
  \STATE $[a^\star_x,T^\star_x] =\argmin_{a\in[0,\bar{\ba}],T\in[\underline{T},\overline{T}]}F(a,T;x,v^\star)$;
  \ENDFOR
 \RETURN $a_x^\star$, $T_x^\star$ for every $x\in\mathbb{Z}$.
 \end{algorithmic}
 \end{algorithm}

Initially, the value function is set to be $v_0(x)=|x-\theta|$. In practice, we only consider the value for a reasonable range of states $x$, say $Z\subset \mathbb{Z}$. We can do this because the transition kernel $q(X(T)=x'|x,a)$ does not have fat tail given a bounded $T$ and goes to $0$ very fast as $|x'-x|$ increases. At iteration $k$, the algorithm does value iteration defined by (\ref{VIFixedRate}) for every $x$ in the reasonable range $Z$. The convergence of the value iteration is guaranteed by the fact that $\min_{a\in[0,\bar{\ba}],T} F(a,T;x,v_k),x\in\mathbb{Z}$ defines a contraction mapping on $v(x),x\in\mathbb{Z}$. Value iterations provide us $v^\star(x)$, i.e., the value function that satisfies the DP equation. Then, we can obtain the optimal action $a^\star_x$ and the optimal waiting time for next observation $T_x^\star$ by computing $[a^\star_x,T^\star_x] =\argmin_{a\in[0,\bar{\ba}],T\in[\underline{T},\overline{T}]}F(a,T;x,v^\star)$.

\begin{remark}
Algorithm \ref{ObConAlgo} can also be applied to solve inhomogeneous arrival process with the finite dimensional optimization in line 3 and 14 replaced by the optimal control problem (\ref{InnerOCExamp}).
\end{remark}

\begin{remark}
Compared with value iterations for classical MDP problems \cite{bertsekas1996neuro} where each iteration, one needs to solve an optimization problem over action variables, the value iterations here requires solving an (infinite-dimensional) optimization problem over both action and observation variables. The convergence rate of value iteartions here is given by $\Vert v_{k+1}-v^\star \Vert\leq \beta^{\underline{T}} \Vert v_{k}-v^\star \Vert$.
\end{remark}
\subsection{Numerical Studies}

In this subsection, we implement Algorithm \ref{ObConAlgo} numerically to show the value of each state $\{v(x)\}_{x\in Z}$, the optimal action $a^\star_x$ given $x$ and the optimal time for the next observation $T^\star_x$.

\begin{figure}[h]
    \centering
    \includegraphics[width=0.485\textwidth]{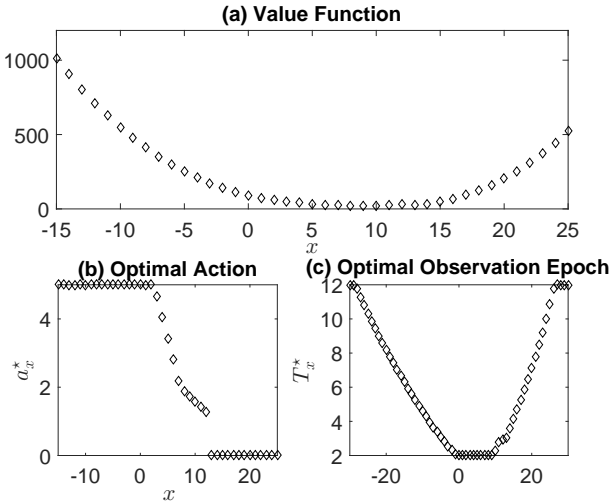}
    \caption{The value function $v^\star(x)$, the optimal action $a^\star_x$ and the optimal time for next observation $T^\star_x$ with respect to $x$.}
    \label{fig:ValueActionObs}
\end{figure}

\begin{figure}[h]
    \centering
    \includegraphics[width=0.489\textwidth]{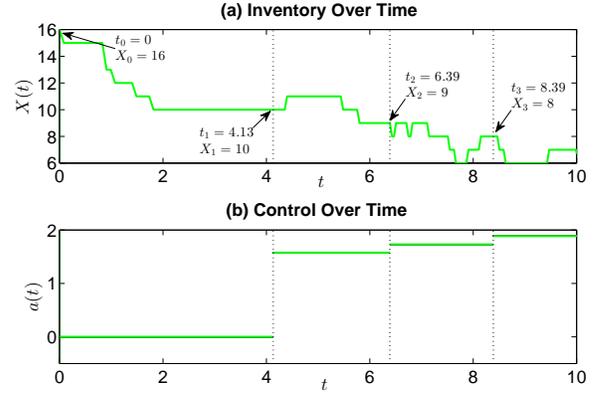}
    \caption{The amount of inventory and the corresponding optimal control and observation over time.}
    \label{fig:PoissonWholeHorizon}
\end{figure}

Here, we set the reference for the amount of inventory to be $\theta = 8$. The departure rate of the Poisson process is $\mu=2$. Here, the Poisson arrival process is homogeneous with upper bound $\bar{\ba}=5$ and lower bound $0$. The time for the next observation is within a range $[\underline{T},\overline{T}]$ where $\underline{T}=2$ and $\overline{T}=12$. The discounted factor is set to be $\beta = 0.8$ and the cost coefficient of accelerating the arrival process $\nu$ is set to be $2$. The cost of observation is set to be $-\kappa T$ where $\kappa=5$. The longer $T$ is, the more frequent the observations are, the lower the cost is. The initial condition when solving the finite-dimensional optimization problem in line $3$ of Algorithm \ref{ObConAlgo} is set to be $a_0=0$ and $T_0=\overline{T}$. 

In Fig. \ref{fig:ValueActionObs}, we present the value, the optimal action, and optimal observation epoch for a range of $x$ near the reference $\theta$. From (a) of Fig. \ref{fig:ValueActionObs}, one should observe that the value function achieves its lowest value near the reference $\theta$. The value function is skewed in a sense that $x$ lower than $\theta$ has higher value than its counterpart who is larger than $\theta$. That is because when $x$ is lower than the reference $\theta$, one has to choose a larger arrival rate $a$ to bring the amount of inventory back to the reference level which would induce more cost. The cost is captured by $a\nu$.

As one can see from (b) of Fig. \ref{fig:ValueActionObs}, when the amount of inventory is low, one chooses the maximum arrival rate to bring the amount of inventory back to the reference level. When the amount of inventory is high, one selects the zero arrival rate to decrease the amount of inventory back to the reference level. At states near the reference, the arrival rate adapts accordingly to keep the amount near the reference while avoiding the acceleration cost.

From (c) of Fig. \ref{fig:ValueActionObs}, we can see the optimal time for next observation given current observation $x$. As we can see, when $x$ near the reference level, the optimal time for next observation is chosen to be $\underline{T}$ which means that one observes more frequently. This is because when $x$ is near $\theta$, the optimal action $a^\star$ adapts constantly as the amount of inventory changes. That means one needs to observe more frequently to get the instant state in order to adapt its action constantly to avoid the acceleration cost or deviation cost. When $x$ is far away from $\theta$, since the optimal action barely adapts, the optimal time for next observation becomes thus longer, which means that one observes less often to avoid observation cost.

Fig. 2 shows how the optimal control and the optimal observation comes into the evolution of the amount of inventory over time. At $t_0=0$ where first observation happens, the decision maker observes that $X_0 = X(t_0) = 16$. From the optimal policy that we obtain from the value iteration, we know $a^\star_{16}=0$ and $T^\star_{16}=4.13$. Thus, the optimal time for the next observation is $t_1= t_0 + T_{16}^\star$ and the optimal Poisson arrival rate for the time interval $[t_0,t_1]$ is $0$. At $t_1$, the decision maker has his second observation which is scheduled at his first observation. The decision maker observes that $X_1 =X(t_1)=10$. By the policy learned before, $a^\star_{10}=1.58$ and $T^\star_{10}=2.26$ which means the next observation time would be $t_2 = t_1 + T^\star_{10}=6.39$ and the optimal Poisson arrival rate is $1.58$ for time interval $[t_1,t_2]$. This procedures go on and on. Since we simulate the homogeneous Poisson arrival rate at each time interval $[t_i,t_{i+1}],i=0,1,\cdots$, the Poisson rate is fixed at each time interval. For inhomogeneous Poisson arrival processes, we have bang-bang type of control at each interval as is shown in (\ref{NecCon}) instead of a fixed one.

\section{Conclusion}
In this paper, we have studied continuous-time jump Markov decision processes with joint control of actions and observation epochs. We have transformed the continuous-time jump MDP model into a regular MDP problem by formulating consolidated utilities between two observation epochs. We have thus obtained the dynamic programming equation with which one can use value iterations to characterize the optimal time for the next observation and the optimal control trajectory. Two case studies have been provided: one is for the gated queueing system where the optimal observation and the optimal action are characterized theoretically. The other one is the inventory control problem with Poisson arrival process in which we have numerically computed the optimal observation epochs and the optimal actions. The results have indicated that the observation is less often at a region of states where the optimal action barely changes. Future works would investigate the finite-horizon jump MDP with a limited number of observations and develop sophisticated learning schemes to learn the optimal action and the optimal observation times when the model is unknown. It is also worth of studying optimal observation for classical MDP model with finite state under a discrete time setting.

\bibliography{references}
\bibliographystyle{IEEEtran}

\end{document}